\documentclass[11pt]{article}
\usepackage{graphicx}
\usepackage{color}
\usepackage{amsmath}
\usepackage{amssymb}
\usepackage{amscd}
\usepackage{framed}
\usepackage{bbm}

\usepackage{fancyhdr,a4wide}
\usepackage{amsthm}

\newcommand{\R}{\mathbb{R}}
\newcommand{\inr}[1]{\left\langle #1 \right\rangle}

\newcommand{\E}{\mathbb{E}}

\newcommand{\PP}{\mathbb{P}}

\newcommand{\eps}{\varepsilon}

\newtheorem{Theorem}{Theorem}[section]
\newtheorem{Lemma}[Theorem]{Lemma}

\newtheorem{Definition}[Theorem]{Definition}

\newtheorem{Corollary}[Theorem]{Corollary}
\newtheorem{Remark}[Theorem]{Remark}

\newtheorem{Question}[Theorem]{Question}

\numberwithin{equation}{section}

\def \proof {\noindent {\bf Proof.}\ \ }

\def \endproof
{{\mbox{}\nolinebreak\hfill\rule{2mm}{2mm}\par\medbreak}}
\def\IND{\mathbbm{1}}

\def\IND{\mathbbm{1}}

\begin{document}

\title{Approximating $L_p$ unit balls via random sampling}
\author{Shahar Mendelson\footnote{Mathematical Sciences Institute, The Australian National University, Canberra, Australia.
\newline
{\sf email: shahar.mendelson@anu.edu.au}
 } }

\medskip
\maketitle
\begin{abstract}
Let $X$ be an isotropic random vector in $\R^d$ that satisfies that for every $v \in S^{d-1}$, $\|\inr{X,v}\|_{L_q} \leq L \|\inr{X,v}\|_{L_p}$ for some $q \geq 2p$. We show that for $0<\eps<1$, a set of $N = c(p,q,\eps) d$ random points, selected independently according to $X$, can be used to construct a $1 \pm \eps$ approximation of the $L_p$ unit ball endowed on $\R^d$ by $X$.  Moreover, $c(p,q,\eps) \leq c^p \eps^{-2}\log(2/\eps)$; when $q=2p$ the approximation is achieved with probability at least $1-2\exp(-cN \eps^2/\log^2(2/\eps))$ and if $q$ is much larger than $p$---say, $q=4p$, the approximation is achieved with probability at least $1-2\exp(-cN \eps^2)$.

In particular, when $X$ is a log-concave random vector, this estimate improves the previous state-of-the-art---that $N=c^\prime(p,\eps) d^{p/2}\log d$ random points are enough, and that the approximation is valid with constant probability.
\end{abstract}

\section{Introduction}
Let $\mu$ be a centred measure on $\R^d$ and for $p \geq 1$ set
$$
B\left(L_p(\mu)\right) = \left\{ v : \in \R^d : \int_{\R^d} |\inr{v,x}|^p d\mu(x) \leq 1 \right\}
$$
to be the unit ball corresponding to the $L_p$ norm endowed on $\R^d$ by $\mu$.

Even at an intuitive level, the sets $B\left(L_p(\mu)\right)$ seem significant because they  ``code" some information on the measure $\mu$. But the fact of the matter is that they are far more important than one might first suspect. Their dual bodies, the so-called $Z_p(\mu)$ bodies, were introduced by E.~Lutwak and G.~Zhang (under different normalization) in \cite{MR1601426}.  G.~Paouris discovered in his seminar work \cite{MR2276533} that when $\mu$ is log-concave\footnote{A log-concave measure has a density that is a log-concave function on $\R^d$.}, the geometry of this family of bodies captures vital information on properties of the generating measure. An alternative, equivalent approach was developed independently by B.~Klartag in \cite{MR2319512}, using the logarithmic Laplace Transform (for a presentation of a unified version of the two approaches, see \cite{MR2852254}).

More information on the geometry of $B\left(L_p(\mu)\right)$ and $Z_p$ bodies and their central role in modern \emph{Asymptotic Geometric Analysis} can be found in the books \cite{MR3331351} and \cite{MR3185453}.
\vskip0.3cm
The identity of the sets $B\left(L_p(\mu)\right)$ is very useful when it comes to the analysis of statistical algorithms involving the measure $\mu$. However, in many statistical applications $\mu$ is not known, and rather than knowing the measure, one is given a sample $X_1,...,X_N$, selected independently according to $\mu$. Thus, it is natural to ask whether the sets $B\left(L_p(\mu)\right)$ can be recovered, or at least approximated using a random sample---hopefully, of a small cardinality, and that the recovery procedure is successful with high probability.

The question of estimating $B\left(L_2(\mu)\right)$ is called \emph{covariance estimation} in statistical literature, and has been studied extensively is recent years (see, e.g. \cite{MR3217437,MR3556768,MR3851758,MR4036049,MR4124338,MenCCM}). Classical results focus on situations where $\mu$ is ``well-behaved", in the sense that linear functionals $\inr{v,\cdot}$ are very light-tailed. Recently, sharp estimates were obtained in heavy-tailed situations. Roughly put, and without going into technical details, the covariance of a centred random vector $X$ can be recovered under very mild assumptions: given a sample $X_1,...,X_N$ for $N \geq c(\eps)d$, one can find ${\cal K}={\cal K}(X_1,...,X_N) \subset \R^d$, such that
$$
(1-\eps){\cal K} \subset B\left(L_2(\mu)\right) \subset (1+\eps){\cal K}.
$$
The linear dependence on $d$ is clearly optimal, while the best estimate that is currently known on $c(\eps)$ is $\sim \eps^{-2}\log(2/\eps)$ (see \cite{MenCCM}). At the same time, all the methods used for covariance estimation are valid only when $p=2$ and do not extend to any other value of $p>2$, even if $X$ is a gaussian random vector, let alone in more general scenarios. As a result, the question of estimating $B\left(L_p(\mu)\right)$ using random data remained completely open, and in the few cases where partial results were known (e.g., Theorem \ref{thm:GR}, below), the estimates were far from satisfactory.
\vskip0.3cm
Since the covariance can be effectively estimated from a sample whose cardinality is proportional to the dimension of the underlying space, let us fix one such  structure:
\begin{Definition}
A measure $\mu$ on $\R^d$ is \emph{isotropic} if it is centred and has the identity as its covariance. In particular, $B(L_2(\mu))=B_2^d$, the Euclidean unit ball. 
\end{Definition}
It is standard to verify that every measure $\mu$ on $\R^d$ has an affine image that is isotropic.
\vskip0.3cm
Given that normalization, the question we wish to address is:
\begin{Question} \label{qu:identification}
Let $\mu$ be an isotropic measure on $\R^d$ and set $X$ to be the random vector distributed according to $\mu$. For $p > 2$ and $0<\eps,\delta<1$ find $N=N(\eps,\delta,p)$ and a mapping $\Phi_p:(\R^d)^N \times \R^d \to \R_+$ for which, with $\mu^N$-probability at least $1-\delta$, for any $v \in \R^d$,
$$
(1-\eps)\E |\inr{X,v}|^p \leq \Phi_p \left( \left(X_i\right)_{i=1}^N,v\right) \leq (1-\eps)\E |\inr{X,v}|^p.
$$
\end{Question}

An example that would be of particular interest is when $\mu$ is an isotropic, log-concave measure, and the state of the art estimate for such measures is due to Gu\'edon and Rudelson:

\begin{Theorem} \label{thm:GR} \cite{MR2304336}
There exists an absolute constant $c$ for which the following holds. Let $\eps \in (0,1)$, $p \geq 2$ and $d \geq d_0(\eps,p)$. Let $X$ be an isotropic, log-concave random vector in $\R^d$ and set
$$
N \geq \frac{(cp)^p}{\eps^2} d^{p/2} \log d.
$$

Then for $t>\eps$, with probability at least $1-2\exp(-c(t/(\eps c_p^\prime))^{1/p})$, for any $v \in \R^d$,
$$
(1-t) \E |\inr{X,v}|^p \leq \frac{1}{N}\sum_{j=1}^N |\inr{X_j,v}|^p \leq (1+t) \E |\inr{X,v}|^p.
$$
\end{Theorem}

Theorem \ref{thm:GR} implies that for log-concave random vectors, $N = c(p,\eps)d^{p/2}\log d$ random points suffice to construct a $1\pm \eps$ approximation that is valid with a constant confidence level---say $\delta=1/2$. Moreover, the approximation procedure is the most natural choice---the $p$-empirical mean,
$$
\frac{1}{N}\sum_{i=1}^N |\inr{X_i,v}|^p.
$$
As a concentration result for $p$-empirical means, the estimate in Theorem \ref{thm:GR} is close to the best that one can hope for. However, the choice of the $p$-empirical mean as an approximation procedure happens to be far from optimal. The problem with that choice can be seen even for a single function: if $f:\Omega \to \R$ is relatively heavy-tailed, a typical sample $(|f(X_i)|^p)_{i=1}^N$ contains enough atypically large values, making the $p$-empirical average too high. As an example, consider $p=2$ and let $f$ be a function for which the Chebychev bound
$$
\PP\left(\left|\frac{1}{N}\sum_{i=1}^N f^2(X_i) - \E f^2 \right| \geq \sqrt{\frac{ \E f^4(X)}{\delta N}} \right) \leq \delta
$$
is sharp. Thus, for the $2$-empirical mean to be a $1\pm \eps$ approximation of $\|f\|_{L_2}^2$ with confidence $1-\delta$, the sample size scales as $1/\delta$. In contrast, as we explain in what follows, the optimal (subgaussian) dependence on $\delta$ can be achieved: the right estimate scales like $\sqrt{\log(2/\delta)}$ (at least in the range we are interested in).
Naturally, the way a recovery procedure performs on a single function says very little on the dependence of $N$ on the dimension in Question \ref{qu:identification}; rather, it indicates that the $p$-empirical mean is likely to be a suboptimal way of approximating the $L_p$ norm. And indeed, our main result is that a different procedure requires a sample size that scales linearly in the dimension rather than like $d^{p/2}\log d$, and performs with very high (subgaussian) probability.

\vskip0.3cm
To formulate the result we require a standard definition.

\begin{Definition}
Let $q > p \geq 1$. The random vector $X$ satisfies an $L_q-L_p$ norm equivalence with constant $L$ if for every $v \in \R^d$, $\|\inr{X,v}\|_{L_{q}} \leq L \|\inr{X,v}\|_{L_{p}}$.
\end{Definition}

\begin{Remark}
One should keep in mind that log-concave random vectors are absolutely continuous and satisfy $L_q-L_p$ norm equivalence with constant $cq/p$ for a suitable absolute constant $c$ (see, for example, \cite{MR3185453}).
\end{Remark}
\vskip0.3cm

The procedure we use is as follows:
\begin{framed}
Given a tuning parameter $\frac{1}{N} \leq \theta <1$ and $p \geq 2$, set
\begin{equation} \label{eq:truncation-p}
\Psi_{p,\theta}(v) =\frac{1}{N} \sum_{j \geq \theta N} \left(|\inr{X_i,v}|^p\right)^*_j
\end{equation}
where $(z_j^*)_{j=1}^N$ is the non-increasing rearrangement of $(|z_j|)_{j=1}^N$.
\end{framed}

\begin{Remark}
Note that $\Psi_{p,\theta}$ is $p$-positive homogeneous; therefore, it suffices to show that
\begin{equation} \label{eq:iso-psi-intro-1}
(1-\eps) \E |\inr{v,X}|^p \leq \Psi_{p,\theta}(v) \leq (1+\eps) \E |\inr{v,X}|^p
\end{equation}
for every $v$ in some Euclidean sphere to ensure that \eqref{eq:iso-psi-intro-1} holds for every $v \in \R^d$.
\end{Remark}

\begin{Theorem} \label{thm:L-p-approx-intro}
There is an absolute constant $c_0$ for which the following holds. Let $X$ be an isotropic random vector in $\R^d$ that is absolutely continuous. For any $\eps \in (0,1)$ and $\theta=c_0\eps^2 N$,
\begin{equation} \label{eq:L-p-approx-intro}
(1-\eps) \E |\inr{X,v}|^p \leq \Psi_{p,\theta}(v) \leq (1+\eps) \E |\inr{X,v}|^p \ \ {\rm \ for \ every \ } v \in \R^d,
\end{equation}
provided that
$$
N \geq c_1^p d \frac{\log(2/\eps)}{\eps^2}
$$
in the following cases:
\begin{description}
\item{$(a)$} If $X$ satisfies an $L_{2p}$-$L_p$ norm equivalence with constant $L$, then $c_1$ depends only on $L$ and \eqref{eq:L-p-approx-intro} holds with probability at least $1-2\exp(-c_2(L,\eps) N)$ for $c_2(L,\eps) = c(L)\frac{\eps^2}{\log^2(2/\eps)}$.
\item{$(b)$} The logarithmic factor in $c_2(L,\eps)$ from $(a)$ is not needed if $X$ satisfies an $L_{q}$-$L_p$ norm equivalence with constant $L$ for some $q>2p$. The claim holds with probability at least $1-2\exp(-c_2^\prime \eps^2 N)$ and $c_2^\prime$ depends on $q-2p$ and $L$.
\item{$(c)$} In particular, if $X$ is log-concave then $c_1$ is an absolute constant and \eqref{eq:L-p-approx-intro} holds with probability at least $1-2\exp(-c_3 \eps^2 N)$ for a suitable absolute constant $c_3$.
\end{description}
\end{Theorem}

Part $(c)$ of Theorem \ref{thm:L-p-approx-intro} follows immediately from Part $(b)$ because a log-concave vector is absolutely continuous and satisfies an $L_{4p}-L_{p}$ norm equivalence with an absolute constant. Thus, the number of vectors that suffice for the construction of a $1 \pm \eps$ approximation of $B(L_p(\mu))$ for a log-concave measure $\mu$ scales linearly in the dimension $d$.

\vskip0.3cm

Observe that for every $\theta$ and $p$, 
$$
\Psi_{p,\theta}(v) \leq \frac{1}{N}\sum_{i=1}^N |\inr{X_i,v}|^p;
$$
thus, Theorem \ref{thm:L-p-approx-intro} leads to a one-sided (lower) bound on the empirical mean:
\begin{Corollary} \label{cor:main}
In the situations described in Theorem \ref{thm:L-p-approx-intro} we have that for every $v \in \R^d$,
$$
(1-\eps) \E |\inr{X,v}|^p \leq \frac{1}{N} \sum_{i=1}^N |\inr{X_i,v}|^p.
$$
\end{Corollary}
Corollary \ref{cor:main} is one in a long line of results which show that one-sided inequalities for the empirical mean are almost universally true---under only minimal assumptions on $X$ (see, for example, \cite{MR3367000,MR3431642,MenCCM} for results of a similar flavour). And while the lower bound is universal, the upper one is highly restrictive, and is false in general. The ``truncation" functional $\Psi$ addresses the problem of atypically large values that are likely to appear in each vector $(|\inr{X_i,v}|)_{i=1}^N$ when $X$ is heavy-tailed, and that leads to the two-sided estimate of Theorem \ref{thm:L-p-approx-intro}. The crucial point is that $\Psi$ endows an adaptive truncation level---based on the nonincreasing rearrangement of the vector $(|\inr{X_i,v}|)_{i=1}^N$---rather than at a fixed value.

\begin{Remark}
The definition of $\Psi_{p,\theta}$ can be extended beyond the class of linear functionals by setting for $f:\Omega \to \R$,
$$
\Psi_{p,\theta}(f) =\frac{1}{N} \sum_{j \geq \theta N} \left( |f(X_i)|^p\right)^*_j.
$$
Moreover, Theorem \ref{thm:L-p-approx-intro} can be extended to far more general function classes than the set of linear functionals on $\R^d$. However, that requires the development of a rather involved technical machinery that is not needed for addressing Question \ref{qu:identification}. We defer the study of the general scenario to \cite{LugMen}, and devote this work to the ``shortest path" leading to the proof of Theorem \ref{thm:L-p-approx-intro}.
\end{Remark}

\vskip0.3cm
In what follows we denote the expectation $\E f(X)$ by $\PP(f)$ and $\PP_N(f)=\frac{1}{N}\sum_{i=1}^N f(X_i)$ is the empirical mean of $f$. We use the same notation---$\PP(A)$ and $\PP_N(A)$ to denote the actual and empirical measures of a set $A$. Absolute constants are denoted by $c$, $C$, etc.; their value may change from line to line. $a \lesssim b$ implies that there is an absolute constant $c$ such that $a \leq c b$. $c(v)$ and $c_v$ denote constants that depend only on the parameter $v$.

\vskip0.3cm

Contrary to what one may expect, the proof of Theorem \ref{thm:L-p-approx-intro} is rather simple. It is based on two facts: first, that for any function $f$, $\Psi_{p,\theta}(f)$ is a sharp estimate of $\E |f|^p$ if $X_1,...,X_N$ satisfies certain ratio estimates of the form
$$
\sup_{\{t : \PP(|f|>t) \geq \eta\}} \left|\frac{\PP_N (|f|>t)}{\PP(|f|>t)}-1\right| < \eps.
$$
Second, that there is an event of high $\mu^N$-probability for which the required ratio estimates are satisfied uniformly by all the linear functionals $F=\{\inr{v,\cdot} : v \in \R^d\}$. This relies heavily on the fact that $F$ is small in an appropriate sense.

The two components of the proof are presented in the next two sections.

\section{Empirical tail integration} \label{sec:sets-to-global}
The goal of this section is to show that, for an arbitrary function $f$, $\Psi_{p,\theta}(f)$ is a good estimator of $\E |f|^p$ for a well-chosen $\theta$, as long as there is enough information on the ratios $\PP_N(|f| \in I)/\PP(|f| \in I)$ for generalized intervals $I$; here an in what follows generalized intervals are open/closed, half-open/closed intervals in $\R$---including rays.

\begin{framed}
Let $0 \leq \lambda, \Delta <1$ and $C \geq 1$. For a function $f$ let ${\cal A}_{\lambda,C,\Delta}$ be the event on which the following holds:
\begin{description}
\item{$(1)$} For any $t>0$ such that $\PP(|f|>t) \geq \Delta$, we have that
$$
\left|\frac{\PP_N (|f| > t)}{\PP (|f|>t)} -1 \right| \leq \lambda;
$$
\item{$(2)$} If $j \in \mathbbm{N}$ and $t>0$ satisfy that $2^{-j}\PP(|f|>t) \geq \Delta$, then
\begin{equation} \label{eq:prop-2}
\left|\frac{\PP_N (|f| > t)}{\PP (|f|>t)} -1 \right| \leq 2^{-j/2};
\end{equation}
\item{$(3)$} For any generalized interval $I \subset \R$,
$$
\PP_N(|f| \in I) \leq \frac{3}{2}\PP(|f| \in I) + C\Delta.
$$
\end{description}
\end{framed}

Ratio estimates are natural in the context of Question \ref{qu:identification} because
$$
\E |f|^p = \int_0^\infty pt^{p-1}\PP(|f|>t)dt.
$$
If sufficiently sharp ratio estimates are available, it is possible to approximate this integral by an empirical functional of the form
$$
\int_0^T pt^{p-1} \PP(|f|>t) dt.
$$

Note that Property $(2)$ is a collection of `isomorphic' estimates that become closer to an isometry for larger sets. For example, if $\PP(|f|>t)$ is of the order of constant, then the allowed distortion in \eqref{eq:prop-2} can be as small as
$$
\left|\frac{\PP_N(|f|>t)}{\PP(|f|>t)}-1\right| \lesssim \sqrt{\Delta},
$$
whereas when $\PP(|f|>t) \sim \Delta$, the allowed distortion in \eqref{eq:prop-2} is $1/2$. That fits the idea of approximating the integral by an empirical counterpart: for a well-chosen $T$, $\int_0^T pt^{p-1} \PP_N(|f|>t)dt$ can be very close to $\int_0^T pt^{p-1} \PP(|f|>t)dt$ even when the distortion is relatively large for sets $\{|f|>t\}$ whose measure is small; at the same time, minimal distortion is essential for sets of relatively large measure, as those have a much higher impact on the two integrals.

\begin{Definition} \label{def:p-error}
For a function $f$, $p \geq 1$ and $T>0$ set
$$
{\cal E}_{T,p}(f)= 2 \sqrt{\Delta}\int_0^T pt^{p-1} \sqrt{\PP(|f|>t)}dt.
$$
Also, for $0<\eta<1$, let
$$
Q_{1-\eta}(f)=\inf\left\{ t : \PP(f>t) < \eta\right\},
$$
i.e., $Q_{1-\eta}(f)$ is the $\eta$ quantile of $f$.
\end{Definition}

\begin{Theorem} \label{thm:main-single-integral}
There are absolute constants $c_1,...,c_5$ for which the following holds. Let $p \geq 1$ and assume that $f(X) \in L_p$ is nonnegative and absolutely continuous. Set $\frac{c_1}{N}< \Delta \leq 1/2$ and assume that $f$ satisfies properties $(1)$-$(3)$ on the event ${\cal A}={\cal A}_{\lambda,C,\Delta}$ for $\lambda=1/2$ and $C=2$. Setting $\theta = c_2 \Delta$ and $\Lambda=Q_{1-c_3\Delta}$ we have that on the event ${\cal A}$,
\begin{equation} \label{eq:upper-single-main}
\frac{1}{N} \sum_{j \geq \theta N} \left(f^p(X_i)\right)_j^* \leq \E f^p + c_4 {\cal E}_{\Lambda,p}(f),
\end{equation}
and
\begin{equation} \label{eq:lower-single-main}
\frac{1}{N} \sum_{j \geq \theta N} \left(f^p(X_i)\right)_j^* \geq \E f^p - c_4 \left({\cal E}_{\Lambda,p}(f)+ \E f^p \IND_{\{f > Q_{1-c_5 \Delta}\}}\right).
\end{equation}
\end{Theorem}

%\begin{Theorem} \label{thm:main-single-integral}
%There are absolute constants $c_1,...,c_6$ for which the following holds. Let $p \geq 1$ and assume that $f(X) \in L_p$ is nonnegative and absolutely continuous. Set $c_1\frac{\log N}{N} \leq \Delta \leq 1/2$, $\theta = c_2 \Delta$ and $\Lambda=Q_{1-c_3\Delta}$. Then with probability at least $1-2\exp(-c_4 \Delta N)$,
%\begin{equation} \label{eq:upper-single-main}
%\frac{1}{N} \sum_{j \geq \theta N} \left(f^p(X_i)\right)_j^* \leq \E f^p + c_5 {\cal E}_{\Lambda,p}(f),
%\end{equation}
%and
%\begin{equation} \label{eq:lower-single-main}
%\frac{1}{N} \sum_{j \geq \theta N} \left(f^p(X_i)\right)_j^* \geq \E f^p - c_5 \left({\cal E}_{\Lambda,p}(f)+ \E |f|^p \IND_{\{f > Q_{1-c_6 \Delta}\}}\right).
%\end{equation}
%\end{Theorem}

\begin{Remark} \label{rem:tail-required}
There is no hope of obtaining an empirical-based estimator of $\E f^p$ if the contribution of the tail $\E f^p\IND_{\{f > Q_{1-\kappa}(f) \}}$ is too big. The reason is that the set $\{f > Q_{1-\kappa}(f)\}$ may be under-represented in the sample: if one is interested in an estimate that holds with $\mu^N$-probability of at least $1-2\exp(-c\Delta N)$, then
$$
\PP_N \left(f > Q_{1-\kappa}(f)\right)=\frac{1}{N}\left|\left\{i: f(X_i) > Q_{1-\kappa}(f) \right\} \right|
$$
may be much smaller than $\kappa=\PP (f > Q_{1-\kappa}(f))$ unless $\kappa$ is of the order of $\Delta$. To see that, let $\kappa \lesssim \Delta$. Then with probability at least $(1-\kappa)^N \geq \exp(-c\kappa N) \geq \exp(-c\Delta N)$ we have that $f(X_i) \leq Q_{1-\kappa}(f)$ for every $1 \leq i \leq N$. On such samples one cannot distinguish between $f$ and of $f \IND_{\{f \leq Q_{1-\kappa}(f)\}}$; however, there can be a significant difference between $\E f^p$ and $\E f^p\IND_{\{f<Q_{1-\kappa}(f)\}}$. As a result, the term $\E f^p \IND_{\{f \geq Q_{1-c_5 \Delta}(f)\}}$ in \eqref{eq:lower-single-main} is essential.
\end{Remark}

Before we turn to the proof of Theorem \ref{thm:main-single-integral}, let us examine the two parameters that are featured in it---namely, ${\cal E}_{T,p}(f)$ for $T=Q_{1-\kappa}(f)$ and $\E f^p \IND_{\{f \geq Q_{1-\kappa}(f)\}}$ for some $0<\kappa<1$. To ease notation we remove the dependence of the two parameters on $f$, and write ${\cal E}_{T,p}$ and $Q_{1-\kappa}$ instead.

\begin{Lemma} \label{lemma:basic-parameters}
Let $f(X) \in L_{2p}$ be nonnegative and absolutely continuous. Then
$$
\E f^p \IND_{\{f > Q_{1-\kappa}\}} \leq \|f\|_{L_{2p}}^p \sqrt{\kappa}
$$
and
$$
{\cal E}_{Q_{1-\kappa},p} \leq c\sqrt{\Delta} \sqrt{\log\left(\frac{1}{\kappa}\right)} \|f\|_{L_{2p}}^p
$$
for an absolute constant $c$.

Moreover, if $f \in L_q$ for $q>2p$ then
$$
{\cal E}_{Q_{1-\kappa},p} \leq c_{q,p}\sqrt{\Delta}\|f\|_{L_{q}}^p
$$
for $c_{q,p}=2p/(q-2p)$.
\end{Lemma}

\proof The proofs of the claims are straightforward. For the first claim observe that
$$
\E f^p \IND_{\{f > Q_{1-\kappa}\}} \leq \left(\E f^{2p}\right)^{1/2} \PP^{1/2}(f > Q_{1-\kappa}) = \|f\|_{L_{2p}}^p \sqrt{\kappa}.
$$
Turning to the two estimates on ${\cal E}_{Q_{1-\kappa},p}$, let $T=Q_{1-\kappa}$ and consider the following two cases. If $T \geq \|f\|_{L_{p}}$ then  by the Cauchy-Schwarz inequality,
\begin{align*}
{\cal E}_{T,p} = & 2\sqrt{\Delta} \left(\int_0^{\|f\|_{L_{p}}} pt^{p-1} \sqrt{\PP(f>t)}dt + \int_{\|f\|_{L_{p}}}^T \frac{1}{\sqrt{t}} \cdot pt^{p-1/2} \sqrt{\PP(|f|>t)}dt \right)
\\
\leq & 2\sqrt{\Delta} \left(\int_0^{\|f\|_{L_{p}}} pt^{p-1}dt + \left(\log\frac{T}{\|f\|_{L_{p}}}\right)^{1/2} \cdot \sqrt{\frac{p}{2}} \left(\int_{\|f\|_{L_{p}}}^T 2p t^{2p-1}\PP(|f|>t)dt\right)^{1/2} \right)
\\
\leq & 2\sqrt{\Delta} \left(\|f\|_{L_p}^p + \left(\frac{p}{2}\log\left(\frac{T}{\|f\|_{L_{p}}}\right)\right)^{1/2}\|f\|_{L_{2p}}^p\right)=(*).
\end{align*}
Recalling that $f$ is absolutely continuous,
$$
\kappa = \PP\left(f \geq Q_{1-\kappa}\right) \leq \left(\frac{\|f\|_{L_p}}{Q_{1-\kappa}}\right)^p.
$$
In particular, $T/\|f\|_{L_p}=Q_{1-\kappa}/\|f\|_{L_p} \leq 1/\kappa^{1/p}$ and
$$
\frac{p}{2}\log\left(\frac{T}{\|f\|_{L_{p}}}\right) \leq \frac{p}{2}\log\left(\frac{1}{\kappa^{1/p}}\right) \leq \frac{1}{2}\log\left(\frac{1}{\kappa}\right).
$$
Therefore,
$$
(*) \lesssim \sqrt{\Delta} \left(\|f\|_{L_p}^p + \|f\|_{L_{2p}}^p \sqrt{\log \left(\frac{1}{\kappa}\right)} \right) \lesssim \sqrt{\Delta \log\left(\frac{1}{\kappa}\right)} \|f\|_{L_{2p}}^p,
$$
as claimed. The proof in the case $Q_{1-\kappa} \leq \|f\|_{L_p}$ requires only the trivial estimate on the integral in $[0,\|f\|_{L_p}]$ used above.

Finally, if $f \in L_q$ for $q>2p$ then
$$
\int_0^{\|f\|_{L_{q}}} pt^{p-1}\sqrt{\PP(f>t)}dt \leq \|f\|_{L_q}^p.
$$
and since $\PP(f>t) \leq (\|f\|_{L_q}/t)^q$, it is evident that
\begin{align*}
\int_{\|f\|_{L_{q}}}^T pt^{p-1}\sqrt{\PP(f>t)}dt \leq & \|f\|_{L_q}^{q/2} \int_{\|f\|_{L_{q}}}^T pt^{p-1-q/2}dt \leq \|f\|_{L_q}^{q/2} \cdot \frac{p}{(q/2)-p} \|f\|_{L_q}^{p-q/2}
\\
= & \frac{2p}{q-2p}\|f\|_{L_q}^p.
\end{align*}
\endproof

Combining Theorem \ref{thm:main-single-integral} with Lemma \ref{lemma:basic-parameters} leads to the following, more user-friendly corollary:
\begin{Corollary} \label{cor:main-single}
There are absolute constants $c_1,...,c_6$ for which the following holds. Set $p \geq 1$ and let $f(X) \in L_{2p}$ be nonnegative and absolutely continuous. Set $\frac{c_1}{N} \leq \Delta \leq 1/2$, let $\theta = c_2 \Delta$, and put
$$
\Psi_{p,\theta}(X_1,...,X_N) = \frac{1}{N} \sum_{j \geq \theta N} \left(|f(X_i)|^p\right)_j^*.
$$
Let ${\cal A}$ be the event on which $f$ satisfies properties $(1)$-$(3)$ with constants $\lambda=1/2$ and $C=2$. Then on the event ${\cal A}$,
\begin{equation} \label{eq:two-sided-single-main-cor}
\left|\Psi_{p,\theta}(X_1,...,X_N)-\E f^p \right| \leq c_4 \sqrt{\Delta \log\left(\frac{1}{\Delta}\right)} \|f\|_{L_{2p}}^p.
\end{equation}
Moreover, if $f \in L_q$ for $q>2p$ then with the same probability,
$$
\left|\Psi_{p,\theta}(X_1,...,X_N) -\E f^p \right| \leq c_{q,p} \sqrt{\Delta}\|f\|_{L_{q}}^p,
$$
where $c_{q,p} \sim p/(q-2p)$.
\end{Corollary}

\begin{framed}
In the context of Theorem \ref{thm:L-p-approx-intro}, we show in what follows that there is a high probability event on which, for \emph{every} $v \in \R^d$, $f_v(X)=|\inr{X,v}|$ satisfies properties $(1)$-$(3)$, as long as $\Delta \geq c\frac{d}{N}\log\left(\frac{eN}{d}\right)$. Once that fact is established (see Section \ref{sec:single}), Theorem \ref{thm:L-p-approx-intro} follows immediately from Corollary \ref{cor:main-single}.
\end{framed}

Let us prove the following version of Theorem \ref{thm:main-single-integral} which gives some freedom in the choice of parameters $\lambda$, $C$ and $\Delta$.

\begin{Theorem} \label{thm:main-integral}
Set $p \geq 1$ and let $f(X) \in L_p$ be nonnegative and absolutely continuous. Assume that $X_1,...,X_N$ is a sample for which $f$ satisfies Properties $(1)-(3)$ with constants $\lambda$, $C$ and $\Delta$. Set $\theta \geq 4\Delta \max\{(1+\lambda),C+3/2\}$ and let
$$
\theta_1=\frac{\theta+2C\Delta}{1-\lambda} \ \ \ {\rm and} \ \ \ \theta_2=\frac{\theta-2C\Delta}{1+\lambda}.
$$
Then
$$
\frac{1}{N}\sum_{j \geq \theta N} (f^p(X_i))_j^* \leq \E f^p + 2\sqrt{\Delta} \int_0^{Q_{1-\theta_2}(f)} pt^{p-1}\sqrt{\PP(f>t)}dt,
$$
and
\begin{align*}
\frac{1}{N}\sum_{j \geq \theta N} (f^p(X_i))_j^* \geq & \E f^p - \left(1+\frac{1}{1-\lambda}\right)\E f^p \IND_{\{f \geq Q_{1-\theta_1}(f)\}}
\\
- & 2\sqrt{\Delta} \int_0^{Q_{1-\theta_1}(f)} 2t \sqrt{\PP(f>t)}dt.
\end{align*}
\end{Theorem}

Theorem \ref{thm:main-single-integral} follows directly from Theorem \ref{thm:main-integral} with the choice of $\lambda=1/2$, $C=2$ and for ${\cal A}_{\lambda,C,\Delta}$ that is the set of samples for which $f(X)$ satisfies Properties $(1)$-$(3)$ with those parameters.

\vskip0.3cm

The proof of Theorem \ref{thm:main-integral} requires several preliminary steps, starting with a straightforward observation: clearly, $\theta_2 \geq 2 \Delta$, and therefore,
\begin{equation} \label{eq:observation}
\PP(f>t) \geq 2\Delta \ \ \ {\rm for \ any \ } 0<t \leq Q_{1-\theta_2}(f).
\end{equation}
As a result, all the level sets $\{f>t\}$ for
$0<t \leq Q_{1-\theta_2}(f)$ satisfy Property $(1)$.

\begin{Lemma} \label{lemma:est-on-hat-Q}
Using the notation of Theorem \ref{thm:main-integral}, let $\hat{Q}=(f(X_i))^*_{\theta N}$. Then
$$
Q_{1-\theta_1}(f) < \hat{Q} < Q_{1-\theta_2}(f).
$$
\end{Lemma}

\proof There are at least $\theta N$ indices $i$ such that $f(X_i) \geq \hat{Q}$; thus
\begin{equation} \label{eq:in-proof-emp-measure-1}
\PP_N (f \geq \hat{Q})=\frac{1}{N} \left|\{ i: f(X_i) \geq \hat{Q}\}\right| \geq \theta.
\end{equation}
Therefore, by property $(3)$ for $I=[\hat{Q},\infty)$,
$$
\theta \leq \PP_N (f \geq \hat{Q}) \leq \frac{3}{2}\PP(f \geq \hat{Q}) + C\Delta,
$$
and as $\theta \geq (C+3/2)\Delta$ it follows that
$$
\PP(f > \hat{Q})=\PP(f \geq \hat{Q}) \geq \frac{2}{3}\left(\theta - C\Delta\right) \geq \Delta.
$$
Hence, by Property $(1)$ for $t= \hat{Q}$,
\begin{equation} \label{eq:in-proof-integral-1}
\left|\frac{\PP_N(f > \hat{Q})}{\PP(f > \hat{Q})}-1 \right| \leq \lambda.
\end{equation}

Next, using Property $(3)$ once again, we have that for any $t>0$ and any $\gamma \leq |t|/2$,
$$
\PP_N\left(f \in [t-\gamma,t]\right) \leq \frac{3}{2} \PP(f \in [t-\gamma,t])+C\Delta.
$$
Taking $\gamma \to 0$ and by the absolute continuity of $f(X)$, $\PP_N(f=t) \leq C\Delta$ for any $t > 0$. In particular, for $t=\hat{Q}$
$$
\left|\{i : f(X_i)=\hat{Q}\}\right| \leq C\Delta N.
$$
Hence,
$$
\theta - 2C\Delta < \PP_N(f > \hat{Q}) < \theta+2C\Delta.
$$
Using \eqref{eq:in-proof-integral-1},
$$
\PP(f > \hat{Q}) \leq \frac{\PP_N(f > \hat{Q})}{1-\lambda} < \frac{\theta+2C\Delta}{1-\lambda}=\theta_1,
$$
and
$$
\PP(f > \hat{Q}) \geq \frac{\PP_N(f > \hat{Q})}{1+\lambda} > \frac{\theta-2C\Delta}{1+\lambda}=\theta_2,
$$
implying that
\begin{equation} \label{eq:est-on-hatQ}
Q_{1-\theta_1} < \hat{Q} < Q_{1-\theta_2},
\end{equation}
as claimed.

\endproof

\begin{Lemma} \label{lemma:int-preliminary}
Let $f(X)$ be nonnegative and absolutely continuous. Using the notation of Theorem \ref{thm:main-integral}, for $p \geq 1$,
\begin{equation} \label{eq:int-preliminary}
\int_0^{Q_{1-\theta_1}} pt^{p-1}\PP_N(f>t)dt - \theta \hat{Q}^p \leq \frac{1}{N}\sum_{j \geq \theta N} \left(f^p(X_i)\right)^*_j \leq \int_0^{Q_{1-\theta_2}} pt^{p-1}\PP_N(f>t)dt
\end{equation}
\end{Lemma}

\proof Recall that $\hat{Q}=\left(f(X_i)\right)^*_{\theta N}$, and therefore,
$$
\frac{1}{N}\sum_{i=1}^N f^p\IND_{\{f \leq \hat{Q}\}}(X_i) - \theta \hat{Q}^p \leq \frac{1}{N}\sum_{j \geq \theta N} (f^p(X_i))_j^* \leq \frac{1}{N}\sum_{i=1}^N f^p\IND_{\{f \leq \hat{Q}\}}(X_i).
$$
By tail integration,
$$
\frac{1}{N}\sum_{i=1}^N f^p\IND_{\{f \leq \hat{Q}\}}(X_i) = \int_0^\infty pt^{p-1} \PP_N \left(f\IND_{\{f \leq \hat{Q}\}} > t\right) dt = \int_0^{\hat{Q}} pt^{p-1}\PP_N (f> t) dt;
$$
Lemma \ref{lemma:est-on-hat-Q} shows that $Q_{1-\theta_1} < \hat{Q} < Q_{1-\theta_2}$ and the wanted estimate follows.
\endproof

To control \eqref{eq:int-preliminary}, let us obtain an estimate on $\int_0^T pt^{p-1} \PP_N(f>t)dt$ that holds as long as the probabilities $\PP(f>t)$, $t \in (0,T)$ are large enough and the nonnegative function $f$ satisfies Property $(2)$. To formulate the claim, recall that
$$
{\cal E}_{T,p}(f)= 2 \sqrt{\Delta}\int_0^T pt^{p-1} \sqrt{\PP(|f|>t)}dt.
$$

\begin{Lemma} \label{lemma:integral-1}
Let $T$ be such that $\PP(f>T) \geq \Delta$. Let $(X_1,...,X_N)$ satisfy that, for any $0<t<T$ and $j \in \mathbbm{N}$,
\begin{equation} \label{eq:iso-in-lemma}
{\rm if} \ \ \ 2^{-j}\PP(f>t) \geq \Delta \ \ {\rm then} \ \ \left|\frac{\PP_N(f>t)}{\PP(f>t)}-1\right| \leq 2^{-j/2}.
\end{equation}
Then
$$
\E f^p\IND_{\{f \leq T\}} - {\cal E}_{T,p}(f) \leq \int_0^T pt^{p-1} \PP_N(f>t)dt \leq \E f^p + {\cal E}_{T,p}(f).
$$
\end{Lemma}

We present the proof for $p=2$ and write ${\cal E}_T$ instead of ${\cal E}_{T,2}(f)$. The proof for $p \not=2$ is identical and is omitted.

\proof For every $T>0$ let $j_T$ be the largest integer such that $\PP(f>T) \geq 2^{j} \Delta $, and set $j_0$ to be the smallest integer such that $2^{j_0} \Delta \geq 1$. Therefore, $2^{j_0-1}\Delta \geq 1/2$ and
$$
\frac{1}{2} \PP(f > T) \leq 2^{j_T}\Delta \leq \PP(f>T).
$$

For $j_T \leq j \leq j_0-1$ let
$$
I_j=\{t > 0 : \Delta 2^{j} \leq \PP(f>t) < \Delta 2^{j+1}\}
$$
and observe that
$$
\bigcup_{j=j_T}^{j_0-1} I_j  \supset (0,T).
$$
Moreover, by \eqref{eq:iso-in-lemma}, for $t \in I_j$
$$
(1-2^{-j/2})\PP(f>t) \leq \PP_N (f> t) \leq (1+2^{-j/2})\PP(f>t),
$$
implying that
\begin{align*}
\int_0^{T} 2t \PP_N (f> t)dt = & \sum_{j=j_T}^{j_0-1} \int_{I_j \cap (0,T)} 2t \PP_N (f> t) dt \leq \sum_{j=j_T}^{j_0-1} \int_{I_j \cap (0,T)} 2t (1+2^{-j/2})\PP (f> t) dt
\\
\leq & \int_0^\infty 2t\PP(f>t)dt + \sum_{j=j_T}^{j_0-1} \int_{I_j \cap (0,T)}2t \cdot \Delta 2^{(j/2)+1} dt
\\
\leq & \E f^2 + 2 \sqrt{\Delta} \sum_{j=j_T}^{j_0-1} \int_{I_j \cap (0,T)}2t \sqrt{\PP(f>t)} dt
\\
\leq & \E f^2 + 2\sqrt{\Delta} \int_0^{T} 2t \sqrt{\PP(f>t)} dt = \E f^2 + {\cal E}_T.
\end{align*}
In the reverse direction and using the same argument,
\begin{equation*}
\int_0^{T} 2t \PP_N (f> t) \geq \int_0^T 2t \PP(f>t)dt - 2\sqrt{\Delta} \int_0^T 2t \sqrt{\PP(f>t)} dt = \E f^2 - {\cal E}_T.
\end{equation*}

\endproof

\noindent{\bf Proof of Theorem \ref{thm:main-integral}.}
Apply Lemma \ref{lemma:integral-1} for $T=Q_{1-\theta_1}(f) \equiv Q_{1-\theta_1}$ and $T=Q_{1-\theta_2}(f) \equiv Q_{1-\theta_2}$, which are both valid choices, as
$$
\PP(f > Q_{1-\theta_1} ) \geq \PP(f > Q_{1-\theta_2} ) \geq \Delta.
$$
Thus,
\begin{equation} \label{eq:emp-upper-bound-1}
\int_0^{Q_{1-\theta_2}} pt^{p-1} \PP_N (f >t) dt \leq \E f^p + {\cal E}_{Q_{1-\theta_2},p},
\end{equation}
and
\begin{align} \label{eq:emp-lower-bound-1}
\int_0^{Q_{1-\theta_1}}  pt^{p-1}\PP_N (f >t) dt & \geq \E f^p \IND_{\{f \leq Q_{1-\theta_1}\}} - {\cal E}_{Q_{1-\theta_1},p}  \nonumber
\\
\geq & \E f^p - \left( \E f^p \IND_{\{f > Q_{1-\theta_1}\}}+{\cal E}_{Q_{1-\theta_1},p}\right).
\end{align}

Next, let us show that
\begin{equation} \label{eq:hat-Q-tail}
\hat{Q}^p \theta \leq \frac{1}{1-\lambda} \E f^p\IND_{\{f \geq Q_{1-\theta_1}\}},
\end{equation}
which, by Lemma \ref{lemma:int-preliminary} completes the proof. To that end, recall that $\PP(f>Q_{1-\theta_2}) \geq \Delta$, and that by Lemma \ref{lemma:est-on-hat-Q},
$$
Q_{1-\theta_1} < \hat{Q} < Q_{1-\theta_2};
$$
thus,
\begin{equation} \label{eq:tail-in-proof-int}
\E f^p \IND_{\{f \geq \hat{Q}\}} \leq \E f^p \IND_{\{f \geq Q_{1-\theta_1}\}}.
\end{equation}
Also, since $\PP_N(f \geq \hat{Q}) \geq \theta$ and $$
\PP(f \geq \hat{Q}) \geq \PP(f \geq Q_{1-\theta_2}) \geq \Delta,
$$
it follows from Property $(1)$ that
$$
\E f^p\IND_{\{f \geq \hat{Q}\}} \geq \hat{Q}^p \PP(f \geq \hat{Q}) \geq \hat{Q}^p (1-\lambda) \PP_N(f \geq \hat{Q}) \geq (1-\lambda) \hat{Q}^p \theta,
$$
proving \eqref{eq:hat-Q-tail}.
\endproof

\section{Ratio estimates for linear functionals} \label{sec:single}
Let us turn to the second component needed for the proof of Theorem \ref{thm:L-p-approx-intro}. Set $F=\left\{ \inr{v,\cdot} : v \in \R^d\right\}$ and denote
\begin{equation} \label{eq:U}
U=\left\{ \IND_{\{ |\inr{v,\cdot}| \in I\}} \ : \ v \in \R^d, \ I \subset \R_+ \ {\rm is \ a \ generalized \ interval } \right\}.
\end{equation}
For a binary-valued function $u$, let $\PP_N(u) = \frac{1}{N} \sum_{i=1}^N u(X_i)$ be its empirical mean and set $\PP(u)= \E u$.

\vskip0.3cm

To complete the proof of Theorem \ref{thm:L-p-approx-intro}, it suffices to find, for $\Delta \geq c_0 \frac{d}{N}\log\left(\frac{eN}{d}\right)$, a high probability event ${\cal A}$ on which every $f \in F$ satisfies Properties $(1)$-$(3)$ for $\lambda=1/2$ and $C=2$. In terms of the class of the indicator functions $U$, we show that:
\begin{framed}
For $\Delta \geq c_0 \frac{d}{N} \log\left(\frac{eN}{d}\right)$, with probability at least $1-2\exp(-c_1\Delta N)$, any $u \in U$ satisfies the following:
\begin{description}
\item{(a)} If $\PP(u) \geq \Delta$ then
$$
\left|\frac{\PP_N (u)}{\PP (u)} -1 \right| \leq \frac{1}{2};
$$
\item{(b)} If $j \in \mathbbm{N}$ and  $2^{-j}\PP(u) \geq \Delta$ then
$$
\left|\frac{\PP_N (u)}{\PP (u)} -1 \right| \leq 2^{-j/2};
$$
\item{(c)} $\PP_N(u) \leq \frac{3}{2}\PP(u) + 2\Delta$.
\end{description}
\end{framed}

The proof of this claim is based on standard tools in empirical processes theory.

\begin{Definition} \label{def:VC}
Let $U$ be a class of $\{0,1\}$-valued functions on $\Omega$. A set $\{x_1,...,x_n\}$ is shattered by $U$ if for every $I \subset \{1,...,n\}$ there is some $u_I \in U$ for which $u_I(x_i)=1$ if $i \in I$ and $u_I(x_i)=0$ otherwise.

The VC dimension of $U$ is the maximal cardinality of a subset of $\Omega$ that is shattered by $U$; it is denoted by ${\rm VC}(U)$.
\end{Definition}
We refer the reader to \cite{MR1385671} for basic facts on VC classes and on the VC-dimension.

\vskip0.3cm
\begin{framed}
The connection between the VC dimension and our problem is that the class of indicator $U$ defined in \eqref{eq:U} satisfies that $VC(U) \leq cd$ for a suitable absolute constant $c$. The proof of this fact can be found, for example, in \cite{MR1385671}.
\end{framed}

The validity of Properties (a)-(c) can be verified for any class of binary valued functions whose VC dimension is at most $d$.

\begin{Theorem} \label{thm:Properties-via-VC}
There are absolute constants $c_0$ and $c_1$ for which the following holds. Let $U$ be a class of binary valued functions, and assume that $VC(U) \leq d$. Then for $c_0 \frac{d}{N} \log\left(\frac{eN}{d}\right) \leq \Delta \leq \frac{1}{2}$, with probability at least $1-2\exp(-c_1 \Delta N)$, the function class $U$ satisfied Properties (a)-(c).
\end{Theorem}

\begin{Remark}
It is very likely that Theorem \ref{thm:Properties-via-VC} is known to experts: ratio estimates of that flavour have been used implicitly in \cite{MR1930004},  and more general ratio estimates, such as Proposition~2.8 in \cite{MR2243881}, can also be used to prove Theorem \ref{thm:Properties-via-VC}. However, we could not locate in literature a simple proof of a suitable version of Theorem \ref{thm:Properties-via-VC}, and the rest of this section is devoted to such a proof.
\end{Remark}

The proof of Theorem \ref{thm:Properties-via-VC} is based on Talagrand's concentration inequality for empirical processes indexed by a class of uniformly bounded functions.

\begin{Theorem} \label{thm:talagrand-conc}
There exist an absolute constant $\kappa$ for which the following holds. Let $U$ be a class of functions, set
$$
\sigma_{U}^2 = \sup_{u \in {U}} {\rm var}(u) \ \ \ {\rm and} \ \ \  b=\sup_{u \in {U}} \|u\|_{L_\infty},
$$
and denote by $(\eps_i)_{i=1}^N$ independent, symmetric, $\{-1,1\}$-valued random variables that are also independent of $(X_i)_{i=1}^N$.

Then for every $x>0$, with probability at least $1-2\exp(-x)$,
$$
\sup_{u \in U} \left|\frac{1}{N}\sum_{i=1}^N u(X_i) - \E u \right| \leq \kappa \left(\E\sup_{u \in U} \left|\frac{1}{N}\sum_{i=1}^N \eps_i u(X_i) \right| + \sigma_{U} \sqrt{\frac{x}{N}}+b\frac{x}{N} \right).
$$

Moreover, if $U$ is a class of binary valued functions and ${\rm VC}(U) \leq d$ then
\begin{equation} \label{eq:VC-estimate}
\E\sup_{u \in U} \left|\frac{1}{N}\sum_{i=1}^N \eps_i uX_i) \right| \leq \kappa\left(\sigma_U \sqrt{\frac{d}{N}\log\left(\frac{e}{\sigma_U}\right)}
+\frac{d}{N}\log\left(\frac{e}{\sigma_U}\right)\right).
\end{equation}
\end{Theorem}
The proof of Theorem \ref{thm:talagrand-conc} can be found in \cite{MR1258865} (see also \cite{MR3185193} for a detailed exposition on related concentration inequalities).

\vskip0.3cm

\noindent {\bf Proof of Theorem \ref{thm:Properties-via-VC}.} Let $c_0$ be a well-chosen absolute constant, set $c_0\frac{d}{N} \log \left(\frac{eN}{d}\right) <\Delta<1$ and let $j \geq 0$ such that $2^j \Delta \leq 1$. Set $\eps_j=2^{-j/2}$ and let
$$
U_j=\{u \in U : 2^j \Delta < \PP(u) \leq 2^{j+1} \Delta\}.
$$
Consider the random variables $\sup_{u \in U_j}  \left|\PP_N(u)-\PP(u) \right|=(*)_j$. For every $u \in U_j$ we have that $\PP(u) \geq 2^{j} \Delta$ and therefore it suffices to show that, with high probability,  $(*)_j \leq \eps_j 2^{j} \Delta = 2^{j/2} \Delta$.

To that end, observe that $\sigma_{U_j}^2 \leq 2^{j+1}\Delta$, set
$$
E_j = \E\sup_{u \in U_j} \left|\frac{1}{N}\sum_{i=1}^N \eps_i u(X_i) \right|,
$$
and clearly $VC(U_j) \leq VC(U) \leq d$. It follows from the second part of Theorem \ref{thm:talagrand-conc} that
$$
E_j \leq \kappa \left(\sqrt{2^{j+1}\Delta\frac{d}{N}\log\left(\frac{e}{2^{j+1} \Delta}\right)}+\frac{d}{N}\log\left(\frac{e}{2^{j+1} \Delta}\right)\right) \leq \frac{1}{4} \eps_j 2^j \Delta = \frac{2^{j/2}\Delta}{4}
$$
by the lower bound on $\Delta$ and the choice of $c_0$.

Let $x_j =c_2\eps^2 2^{j}\Delta N=c_2\Delta N$. Invoking the first part of Theorem \ref{thm:talagrand-conc}, it is evident that with probability at least $1-2\exp(-c_3 \Delta N)$
$$
\sup_{u \in U_j}  \left|\PP_N(u)-\PP(u) \right| \leq 2^{-j/2}\Delta,
$$
as required. Thus, Property (b) follows with the wanted probability thanks to union bound for $\{ j \geq 0 : 2^j \Delta \leq 1\}$ and recalling that $\Delta \geq c_0\frac{d}{N} \log \left(\frac{eN}{d}\right)$ for a well-chosen absolute constant $c_0$.

Repeating the same argument for $\eps_j=\lambda$ and $x_j = c_4 \lambda^2 2^j \Delta$ and using the union bound once again, it is evident that with probability at least $1-2\exp(-c_5 \lambda^2 \Delta N)$,
$$
\sup_{\{u \in U : \PP(u) \geq \Delta\}} \left|\frac{\PP_N(u)}{\PP(u)}-1 \right| \leq \lambda.
$$
Property (a) is verified with the wanted probability by setting $\lambda=1/2$.

Finally, let us turn to Property (c). When $\PP(u) \geq \Delta$, Property (c) follows from Property (a), and when $\PP(u) \leq \Delta$ we use Talagrand's concentration inequality again. Indeed, by the lower bound on $\Delta$,
$$
\E\sup_{\{u \in U : \PP(u) \leq \Delta\}} \left|\frac{1}{N}\sum_{i=1}^N \eps_i u(X_i) \right| \leq \kappa \left( \sqrt{\Delta\frac{d}{N}\log\left(\frac{e}{\Delta}\right)}+\frac{d}{N}\log\left(\frac{e}{ \Delta}\right)\right) \leq \frac{\Delta}{2\kappa}
$$
Therefore, with probability at least $1-\exp(-x)$,
$$
\sup_{\{u \in U : \PP(u) \leq \Delta\}} \left|\PP_N(u)-\PP(u)\right| \leq \kappa \left(\frac{\Delta}{2\kappa} + \sqrt{\frac{\Delta x}{N}} + \frac{x}{N}\right) \leq 2\Delta
$$
by setting $x=c_5 \Delta N$. Hence, Property (c) is verified with $C=2$ and with the wanted probability.
\endproof

\vskip0.3cm

\noindent{\bf Proof of Theorem \ref{thm:L-p-approx-intro}.} With probability at least $1-2\exp(-c\Delta N)$, every indicator function in $U$ satisfies Properties (a)-(c). As a result, on that event, for every $v \in \R^d$, $|\inr{v,X}|$ satisfies Properties $(1)$-$(3)$. Now Theorem \ref{thm:L-p-approx-intro} follows immediately from Corollary \ref{cor:main-single}.
\endproof

\section{Concluding Remarks}
The proof of Theorem \ref{thm:L-p-approx-intro} is relatively straightforward, but that is due to good fortune---that the class of indicator functions
$$
U=\left\{ \IND_{\{ |\inr{v,\cdot}| \in I\}} \ : \ v \in \R^d, \ I \subset \R_+ \ {\rm is \ a \ generalized \ interval } \right\}.
$$
is very simple --- it has VC-dimension that is proportional to the algebraic dimension of the underlying space. For more general classes of functions the situation is far more complex: the class of indicators generated by tails of functions in $F$ need not have a finite VC dimension, let alone a well behaved one.  In \cite{LugMen} we develop a theory that allows one to overcome that obstacle. We show that under minimal condition on the class $F$ and with high probability, Properties $(1)$-$(3)$ hold uniformly in the class. As a result, a more general version of Theorem \ref{thm:L-p-approx-intro} happens to be true.

\bibliographystyle{plain}
\bibliography{L-p}

\end{document}